\documentclass[10pt]{article}

\usepackage{amsmath,latexsym,amssymb,amsthm}

\usepackage{tikz}

\numberwithin{equation}{section} \allowdisplaybreaks

\newtheorem{theorem}{\bf\normalsize Theorem}[section]
\newtheorem{corollary}[theorem]{\bf\normalsize Corollary}

\newtheorem{lemma}[theorem]{\bf\normalsize Lemma}
\newtheorem{proposition}[theorem]{\bf\normalsize Proposition}
\newtheorem{remark}[theorem]{\bf\normalsize Remark}

\begin{document}

\title{On tail probability of the  covariance matrix in Eldan's stochastic localization  }
\author{QingYang Guan\footnote{
 Institute of Applied Mathematics,  AMSS, CAS.
   }
   }
\date{} \maketitle

\begin{abstract}

The   Eldan's stochastic localization is a new kind of stochastic evolution  in the space of probability measures which
provides a novel  way to study high dimensional convex body.
A central object in the study of the stochastic localization is the stochastic  process of its covariance matrix.
The main result of  this paper is    some  exponential-type tail probability estimate of the covariance process   for the general time.
This estimate implies a weaker version of a $p$-moment conjecture by Klartag and Lehec. The stochastic localization considered here is a simplified  version by Lee and Vemplala.
\end{abstract}


{\textbf{\small\noindent AMS(2000) Subject Classification:
  Primary $52\mathrm{A}23$;\  Secondary
  $60\mathrm{H}30$ }

\section{Introduction}

Let $\mu$ be a log-concave probability measure with    density  $\rho$ on   $\mathbb{R}^n$, $n\geq 1$. It  is called isotropic
if its barycenter is zero  and its  covariance matrix is the identity, respectively. For $t>0$, $\mu$  is called $t$-uniformly log-concave if
  $\rho(x)=e^{-\varphi(x)-t|x|^2/2}$ for some $\mathbb{R}\cup \{+\infty\}$ valued convex function $\varphi$.

The   Eldan's stochastic localization    of  $\mu$ is a probability distribution  valued stochastic process  $(\mu_t)_{t\geq 0}$ taking
$(p_t)_{t\geq 0}$
as the density process  with   $\mu_0=\mu$ and $p_0=\rho$.
In
a simplified  version by Lee and Vemplala \cite{LV24}, it is the solution to   the following stochastic differential  equation with the initial condition $p_0$:
\begin{align}
  dp_t(x)= p_t(x) \langle x-a_t,d B_t\rangle ,\ \ \ \ \ \   x\in \mathbb{R}^n,\label{p}
 \end{align}
where   $a_t=  \int
   x p_t(x)dx $ is the barycenter of $p_t$ and $(B_t)$ is  a standard Brownian motion in $\mathbb{R}^n$, respectively.
Denote  the covariance matrix  of $p_t$ by
\begin{align}
    A_{t}=&\int  (x-a_t)\otimes (x-a_t)p_t(x) dx,\ \ \ \ t\geq 0.\nonumber
 \end{align}
 For a   complete  statement of this model and its  previous development,
we refer to   Klartag and  Lehec \cite{KL24} and the references therein.
For some  most recent applications of stochastic localization on log-concave measures, we refer to  \cite{G24},  Klartag and Lehec \cite{KL25}\cite{KL25B} and Bizeul \cite{B25}.

It is known   from the previous works that $(A_t)$ is a  central object in the study of the stochastic localization.
 In May 2025, professors  Klartag and    Lehec  formulated the  following question to the author: ``whether the method in \cite{G24}
 can be extended to yield a proof of the following estimate"
\begin{align}
 \mathbb{E}\mathrm{Tr}( A_t^{  p})
\leq (Cp)^{\beta p}n,\ \ \ \ t\geq 0, \ p\geq 1,\ \beta=1. \label{pm}
\end{align}
The author  was told    that this  inequality   held  for a universal constant $\beta>0$, which
 can be verified   by a modification of the method in \cite{G24} (the case for $p=2$).

With many efforts, the author did  not find a way to  prove this $p$-moment conjecture. However,   a new way to look at it shows
 that (\ref{pm}) holds for $\beta=2$.  The content  of this  paper is a combination of
  two drafts on this topic written in  May 2025  to professors  Klartag and Lehec   for their   interests in   $p$-moment estimate. The first draft given in the Appendix below   is to show that
   a modification of the method in \cite{G24}
  may not be a good choice. The second draft provides  the main context of this paper with the same   title.

Instead of considering the $p$-moment problem,   we
study    some  tail probability estimate related to  the covariance process  $(A_t)$.
More explicitly, we consider  the following  problem   which is  stronger   than (\ref{pm}):   whether there exits  some universal constant $c>0$ such that for any isotropic  log-concave probability measure $\mu$ on   $\mathbb{R}^n$,
  \begin{align}
   \sum_{i=1}^n\mathbb{P}  \big(\lambda_{i;t}> 2 \big)
 \leq &   2 e^{-ct^{-1}} n,\ \ \ \ \ t>0,\label{bba1}
\end{align}
where  $0\leq \lambda_{1;t} \leq  \cdots \leq  \lambda_{n;t}$
 are the eigenvalues of the covariance matrix $A_t$, repeated according to their multiplicity. Notice that (\ref{bba1}) is a stronger  form  of     Corollary 3.1 in  \cite{G24}

The following result  is a   weaker version of (\ref{bba1}) and  is also a stronger version of Corollary 3.1 in  \cite{G24}.
Here we do not get an  estimate
for the general tails.
\begin{proposition}\label{Step0}
 There exists a universal constant $c>0$ such that for any isotropic  log-concave probability measure $\mu$ on   $\mathbb{R}^n$ with $n\geq1$,
  \begin{align}
   \sum_{i=1}^n\mathbb{P}  \big(\lambda_{i;t}> \frac{8}{3} \big)
 \leq &   2 e^{-ct^{-1/2}} n,\ \ \ \ \ t>0.\label{ba1}
\end{align}

\end{proposition}

The tail probability estimate of $(A_t)$ began  in Eldan \cite{E13} for the original model.
The following   estimate  of     $(A_t)$ is given   in     Lee and  Vempala \cite{LV24} for the simplified model in (\ref{p}):
\begin{align} \mathbb{P}(\max_{0\leq s\leq t}\|A_s-\mathrm{Id}\|_{op}\geq \frac{1}{2})\leq & 4e^{-(Ct)^{-1}},\ \ \ \ t\in (0,(2Cn^{1/2})^{-1}].\label{LV}  \end{align}
The best known estimate of $\|A_t\|_{op}$  below can be found in Theorem 60 of Klartag and Lehec \cite{KL24},
\begin{align}
 \mathbb{P}(\|A_t\|_{op}\geq 2)\leq e^{-(C_1t)^{-1}}, \ \ \ \ \ t\in[0, C_1^{-1}(\log n)^{-2}] ,\ \ \ \  n\geq 2,\label{oper}
\end{align}
where $C_1>0$ is some universal constant.  See also  Chen  \cite{C21}, Klartag and Lehec \cite{KL22}  and  Bizeul \cite{B24} for more related results.

The inequality in (\ref{ba1})  can be taken  as a generalization of the previous operator norm estimates of $(A_t)$
introduced  above.
Up to now we do not know whether (\ref{bba1}) holds    for $t\in [(\log n)^{-2},(\log n)^{-1}]$,
 the validity of which would improve  the best known isoperimetric    bound $C\sqrt{\log n}$ given in Klartag \cite{K23}.
See  Klartag and  Lehec  \cite{KL24} for more details.
Moreover, (\ref{bba1}) is true under    the KLS
  conjecture  formulated  in   Kannan,
Lov$\acute{\mathrm{a}}$sz, and Simonovits \cite{KLS95}. See Remark \ref{r}  below for more details.

Applying  Proposition \ref{Step0},  we show that (\ref{pm}) holds for $\beta=2$ in Corollary \ref{Step000} below. The advantage of this way to estimate
the $p$-moment    is that we do not need to consider  the   higher power potentials.
The proof of Proposition \ref{Step0} follows the arguments in \cite{G24} with more careful analysis.

Throughout the paper,
the  letters $C,c$ denote positive universal constants, whose value may change from one line to the next. Notation $\log$
stands for the natural logarithm.      For an $n$ by $n$ matrix $A$,   $\|A\|_{op}$ stands for  the   operator norm. We also denote the trace of $A$ by $\mathrm{Tr} (A)$ and
 $\mathrm{Id}$   for the identity matrix, respectively. For a function $f$ defined on the real line,   $f''$ stands for
its   second derivative.    Notation $\mathbb{E}$ refers to the expectation operation and $\int$ is always for $\int_{\mathbb{R}^n}$. For a  set $H$,
the  indicator function of it  {is denoted by $I_H$}.

\section{Proofs}

We use the method in the proof of Lemma 2.1 in \cite{G24} to prove Proposition \ref{Step0}. The following two Lemmas are taken from Lemma 2.2 and (2.23)   in \cite{G24}, respectively.

\begin{lemma}\label{Step6}
Let   $D_0>4$ and    $r_0\in[ \frac{7}{3},\frac{8}{3}]$.  Then there exists a constant $b$ and an
increasing function $f \in C^2(\mathbb{R})$ such that
   \begin{align}\label{fr}
   f(r)=& e^{D_0(r-r_0)},\ \ \ \ \ r\leq r_0-D_0^{-1};\\ \ \ \ \ f(r)=&br^2,\ \ \ \ \ \ \ \ \ \ \ \ r\geq r_0, \label{fr2}
   \end{align}
   and
     \begin{align}
  \frac{1}{20}\leq b\leq \frac{1}{5};\ \ \ \ \ \ \  |f''(r)|\leq   D_0^2f(r) ,\ \ \ \ r\in \mathbb{R}. \label{sece}
   \end{align}
\end{lemma}

 \begin{lemma}\label{Step4}
Let   $\mu$ be an isotropic  log-concave probability measure on   $\mathbb{R}^n$.
Let $D_0>4 $,    $r_0\in [ \frac{7}{3},\frac{8}{3}]$ and define a function $f$   by  Lemma \ref{Step6}.
Set
\begin{align}
F_t=\sum_{i=1}^n f(\lambda_{i;t})  ,\ \ \ \ t\geq 0,\label{ft}
\end{align}
where  $0\leq \lambda_{1;t} \leq  \cdots \leq  \lambda_{n;t}$
 are the eigenvalues of the covariance matrix $A_t$, repeated according to their multiplicity.
   Then, for $t_0\in (0,1]$,
  \begin{align}
 \frac{d}{dt}  \mathbb{E} F_t\leq \big(500t^{-1}+144\sqrt{3}D_0^2 t^{-1/2}\big) \mathbb{E} F_t,\ \ \ \ \ t\in [t_0,1]. \label{ftf}\end{align}
\end{lemma}
 When $\mu$ is  $t$-uniformly log-concave for some $t>0$, we have by
the    log-concave Lichnerowicz inequality (see, e.g., \cite{BGL14})
\begin{align}\label{Lic}
A_t\leq \frac{1}{t} ,\ \ \ \ \ t>0.
 \end{align}  The inequality (\ref{Lic}) is improved by the following
 log-concave spectral-variance  inequality
   in      \cite{K23},
\begin{align}\label{kla}
 C_p(\mu)\leq \sqrt{ \frac{\|Cov(\mu)\|_{op} }{t}} ,\ \ \ \ \ t>0.
 \end{align}

Next we prove  Proposition \ref{Step0}.
\begin{proof}

By considering   some   product of     $\mu$
in a higher dimensional space, we only need to verify (\ref{ba1}) for $n$ big enough.
Therefore  it suffices to  show that  for some universal constant $c\in (0,1)$ the following estimate
holds:
  \begin{align}
  \mathbb{E} \sum_{ \lambda_{i;t}> \frac{8}{3}} \lambda_{i;t}^2
 \leq &  e^{-ct^{-1/2}}  n,\ \ \ \ \ t\in[0,c],\ \ n\geq c^{-1}.\label{tp}
\end{align}
In what below  $\sum_{ \lambda_{i;t}> \frac{8}{3}}$, for example, represents   $\sum_{i:1\leq i\leq n, \lambda_{i;t}> \frac{8}{3}}$.  This convention
  will be used in the rest part of the paper.

 With $p$ defined by   $t=10^{-8}p^{-2}$, (\ref{tp}) is again  a consequence of the following inequality
\begin{align}
  \mathbb{E} \sum_{ \lambda_{i; 10^{-8}p^{-2}}> \frac{8}{3}} \lambda_{i;10^{-8}p^{-2}}^2
 \leq &  e^{-p/20}  n,\ \ \ \ \ \ p>C,\ n>C, \label{zba1}
\end{align}
 where $C>0$ is a universal constant. The rest part of the proof is to verify (\ref{zba1}) and the arguments are divided into four parts.
 We will first  prove some pre-estimate at time $p^{-8}$ in (\ref{k0}) and then apply it to prove   (\ref{zba1}).
 In   part 2), this pre-estimate is known  because the value of  $p$ is  large. \medskip

$1)\ \mathbf{a\ known\ bound\ for \  the \ operator\ norm}$\ \

The estimate in (\ref{oper}) will be used in what below. We remark that the estimate in (\ref{oper}) can also be replaced by (\ref{LV}) for our purpose.   Set  \begin{align}
t^*= C_1^{-1}(\log n)^{-2} ,\ \ \ \  n\geq 2,\nonumber
\end{align}
where  $C_1$ is given in (\ref{oper}).

In what below we assume that $n$ and $p$ are both bigger than a universal constant $C>e^{125}$ such that
\begin{align}
 e^{ 10^4 \log p + 11p/18 }\big( e^{-2p /3}+   e^{-(C_1t_1^*)^{-1} }\big)\leq&  e^{-p /20},\ \ \  \label{p1s}\\
 (C_1t)^{-1}-(10^3+2)|\log t|\geq& |\log t|^3  ,\ \ \ \ \ t\in (0,t^*], \label{nc}\end{align}
where
\begin{align}
 t_1^*=p^{- 8}.
   \end{align}
The assumption $p\geq e^{125}$ shows that
   \begin{align}
 t_1^*\leq e^{-10^3}.\label{t12}
   \end{align}

$2)\ \mathbf{estimates\ for\ the\ expectation\  in} $  (\ref{zba1}) $  \mathbf{when}$ $t_1^* \leq t^*$\ \

In this subsection we assume that $t_1^* \leq t^*$.
 Set
   \begin{align}
  F_t= \sum_{i=1}^nf(\lambda_{i;t}) ,\ \ \ \ \ t\geq 0,\label{ap}
   \end{align}
   where $f$ is a function     in Lemma \ref{Step6} with $r_0= \frac{8}{3}$ and $D_0=p$. \medskip

Notice that $10^{-8}p^{-2}\geq t_1^*$ since $p>e^{125}$. By   Lemma \ref{Step6}, we see that $f$ is an increasing function and $f(2)\leq e^{-2p/3}$.  Then we have by Lemma
\ref{Step4} with   $t_0=t_1^*$, (\ref{Lic}) and  (\ref{oper})
  \begin{align}
    \mathbb{E} F_{10^{-8}p^{-2}}  \leq   & e^{-500 \log t_1^* + p/9 }   \mathbb{E} F_{t_1^*} \label{at} \\
    = & e^{4\cdot 10^3 \log p + p/9 }\mathbb{E} F_{t_1^*} I_{\|A_{t_1^*}\|_{op}\leq 2}+  e^{4\cdot 10^3 \log p  + p/9 } \mathbb{E} F_{t_1^*} I_{\|A_{t_1^*}\|_{op}> 2}\nonumber\\
\leq &e^{10^4 \log p + p/9 }\sum_{i=1}^n  \mathbb{E}  f(\lambda_{i;t_1^*}) I_{\lambda_{i;t_1^*}\leq 2 }+e^{10^4 \log p + p/9 }  \mathbb{P}(\|A_{t_1^*}\|_{op}\geq 2) n \nonumber\\
\leq &e^{10^4 \log p + p/9 }   e^{-2p /3}n+ e^{10^4 \log p + p/9 } e^{-(C_1t_1^*)^{-1} }n\nonumber\\
 \leq  &    20^{-1} e^{-p/20 }   n,\nonumber
   \end{align}
   where   (\ref{p1s}) is  used  in the last  inequality above.
   The estimate above and the lower bound of $b$ in (\ref{sece}) show that  (\ref{zba1}) holds when $t_1^* \leq t^*$.

\medskip

$3)\ \mathbf{  assistant\ stochastic\ processes }$ $  \mathbf{when}$ $t_1^*> t^*$\ \

In the rest part of the proof we always assume that $t_1^*> t^*$. To prove a pre-estimate   at time $t_1^*=p^{-8}$, we first define
a backward sequence $(t_k^*)$.
Set \begin{align}
\  \ \ \   t_{k+1}^*=e^{-(t_k^*)^{-1/16}},\ \ \ \ \ \ k\geq 1.\label{ttkk}
   \end{align}
  Define
      \begin{align}
  k_0=\sup \{k\geq 1: t_k\leq t^*\},\nonumber
   \end{align}
and
      \begin{align}
  t_k=t^*_{k_0-k+1},\ \ \ \  k=1,\cdots,k_0.\nonumber
   \end{align}
   We have
      \begin{align}
\  \ \ \   t_{k+1}=|\log t_{k}|^{-16},\ \ \ \ \ \ k=1,\cdots, k_0-1.\label{ttkk+1}
   \end{align}
From the definitions above we see that
   \begin{align}
 t_1\leq t^*,\ \ \ \ t_{k_0}=t_1^*.\label{ttkk+12}
   \end{align}

Set
   \begin{align}
  s_1=\frac{7}{3}, \ \ \   s_{k+1}= \frac{7}{3}+ \sum_{i=2}^{k+1}|\log t_{i}|^{-1/2}  ,\ \ \ \ \ \ k= 1,\cdots,k_0-1.\nonumber
   \end{align}
With (\ref{t12}), we can check that
      \begin{align}
  \frac{7}{3} \leq  s_1\leq s_{k_0-1} \leq
 \frac{15}{6}.\label{inc}
   \end{align}
This allows us to  define
   \begin{align}
  F_{k;t}= \sum_{i=1}^nf_{k}(\lambda_{i;t}) ,\ \ \ \ \ t\geq 0,\ 1\leq k\leq k_0-1,\label{1ap}
   \end{align}
   where $f_k$ is a function     in Lemma \ref{Step6} with $r_0=s_k$ and $D_0=|\log t_k|^4$.
 \medskip

$4)\ \mathbf{estimates\ for\ the\ expectation\  in} $  (\ref{zba1})  $  \mathbf{ by \ finite \ induction}$ $  \mathbf{when}$ $t_1^*> t^*$\ \ \medskip

$4.1)\  {estimates\ for\ the\ expectation\  of} $   $(F_{1;t})$

By definition, $144\sqrt{3}D_0^2 t^{-1/2} \leq 500t^{-1}$ when $t\in [t_1,t_2]$ and $D_0=|\log t_1|^4$.  For $t\in [t_1,t_2]$, we have by Lemma
\ref{Step4} with $f=f_1$, $t_0=t_1$, Lemma \ref{Step6}, (\ref{oper}) and (\ref{nc})
  \begin{align}
    \mathbb{E} F_{1,t}  \leq   & t_1^{-10^3}  \mathbb{E} F_{1,t_1} \nonumber\\
    = & t_1^{-10^3} \mathbb{E} F_{1,t_1} I_{\|A_{t_1}\|_{op}\leq 2}+   t_1^{-10^3} \mathbb{E} F_{1,t_1} I_{\|A_{t_1}\|_{op}> 2}\nonumber\\
\leq &t_1^{-10^3} \sum_{i=1}^n  \mathbb{E}  f_1(\lambda_{i;t_1}) I_{\lambda_{i;t_1}\leq 2 }+t_1^{-10^3-2}     \mathbb{P}(\|A_{t_1}\|_{op}\geq 2)n\nonumber\\
\leq &t_1^{-10^3}     f_1(2)n+  e^{-(C_1t_1)^{-1}+(10^3+2)|\log t_1|}n\nonumber\\
    = &  e^{-\frac{1}{3}|\log t_1|^4+10^3 |\log t_1| }n+  e^{- |\log t_1|^3  }n \nonumber\\
    \leq   &  e^{-|\log t_1|^2}n,\label{0}
   \end{align}
where $t_1\leq t_1^*$ and (\ref{t12}) are used in the last step above.\medskip

   $4.2)\  {estimates\ for\ the\ expectation\  of} $   $(F_{k;t})$ for $2\leq k\leq k_0-1$

 We claim further that      for every integer $k\in \{1,\cdots,k_0-1\}$,
     \begin{align}
       \mathbb{E} F_{k,t}
    \leq &  e^{-|\log t_{k}|^2}n, \ \ \ \ \ \ t\in [t_{k},t_{k+1}].\label{intr}
   \end{align}  By (\ref{0}),
   the estimates above hold when $k=1$. The rest part of the  proof is the same as
  the induction steps in the proof of Lemma 2.1 in \cite{G24}. \medskip

    $4.3)\  {estimates\ for\ the\ expectation\  of} $   $F_{k_0;10^{-8}p^{-2}}$

  We have by (\ref{intr})  with (\ref{ttkk+1}) and the equality in (\ref{ttkk+12})
     \begin{align}
       \mathbb{E} F_{k_0-1,t_{k_0}}
    \leq    e^{-|\log t_{k_0-1}|^2}n= e^{-(t_1^*)^{-1/8}}n=  e^{-p}n. \label{k0}
   \end{align}

Set
   \begin{align}
  F_{k_0;t}= \sum_{i=1}^nf_{k_0}(\lambda_{i;t}) ,\ \ \ \ \ t\geq 0,\label{t0}
   \end{align}
   where $f_{k_0}$ is a function     in Lemma \ref{Step6} with $r_0=8/3$ and $D_0=p$.
Denote the corresponding  parameter  $b$ of $f_k$ in Lemma \ref{Step6}  by $b_k$ in what below.

We have by Lemma  \ref{Step6} and  (\ref{inc})
\begin{align}
   &   \sum_{i=1}^n  \mathbb{E} f_{k_0}(\lambda_{i;t_{k_0}}) I_{\lambda_{i;t_{k_0}}> s_{k_0-1}} \nonumber\\
     = &\sum_{i=1}^n  \mathbb{E} f_{k_0}(\lambda_{i;t_{k_0}}) I_{s_{k_0-1}<\lambda_{i;t_{k_0}}\leq  8/3}+\sum_{i=1}^n  \mathbb{E} f_{k_0}(\lambda_{i;t_{k_0}}) I_{\lambda_{i;t_{k_0}}> 8/3}\nonumber\\
  \leq &  \sum_{i=1}^n  \mathbb{E}b_{k_0} (\frac{8}{3})^2 I_{s_{k_0-1}<\lambda_{i;t_{k_0}}\leq  8/3}+
      \sum_{i=1}^n \mathbb{E}  b_{k_0}\lambda_{i;t_{k_0}}^2I_{\lambda_{i;t_{k_0}}> 8/3}\nonumber\\
      \leq & 6 \sum_{i=1}^n  \mathbb{E}b_{k_0-1} \lambda_{i;t_{k_0}}^2 I_{s_{k_0-1}<\lambda_{i;t_{k_0}}\leq  8/3}+
      4\sum_{i=1}^n \mathbb{E}  b_{k_0-1}\lambda_{i;t_{k_0}}^2I_{\lambda_{i;t_{k_0}}> 8/3}\nonumber\\  \leq& 6 \mathbb{E} F_{k_0-1,t_{k_0} },\label{ab}
   \end{align}
   where the   bound of $b$ in (\ref{sece})  is used in the second inequality above.
 It follows by    Lemma \ref{Step4} with $f=f_{k_0}$,
$t_0=t_{k_0}$, the definition of $f_{k_0}$ and (\ref{ab}) that
 \begin{align}
   & \mathbb{E} F_{k_0,10^{-8}p^{-2}}\nonumber\\\leq &
    e^{-500 \log t_1^* + p/9 }   \mathbb{E} F_{k_0,t_{k_0}} \label{at2}\\
      =& e^{4\cdot 10^3 \log p + p/9 } \sum_{i=1}^n  \mathbb{E} f_{k_0}(\lambda_{i;t_{k_0}}) I_{\lambda_{i;t_{k_0}}\leq s_{k_0-1}}+ e^{4\cdot 10^3 \log p + p/9 }\sum_{i=1}^n
       \mathbb{E}  f_{k_0}(\lambda_{i;t_{k_0}}) I_{\lambda_{i;t_{k_0}}> s_{k_0-1}}\nonumber\\
     \leq &e^{4\cdot 10^3 \log p + p/9 } f_{k_0}(  s_{k_0-1})n  +e^{5\cdot 10^3 \log p + p/9 } \mathbb{E} F_{k_0-1,t_{k_0} }\nonumber  .\label{k+1}
   \end{align}
  We see that $f_{k_0}(  s_{k_0-1})\leq e^{-p/6}$ by the last inequality in (\ref{inc}) and the definition of $f_{t_0}$ below  (\ref{t0}).  Then the    estimate above, (\ref{p1s}) and (\ref{k0}) show that  \begin{align}
    \mathbb{E} F_{k_0,10^{-8}p^{-2}}
     \leq&   e^{4\cdot 10^3 \log p + p/9 }  e^{-p/6}n+ e^{5\cdot 10^3 \log p + p/9 }  e^{-p}n \nonumber\\
    \leq  &       20^{-1}e^{-p/20}n ,\nonumber
   \end{align}
   which implies (\ref{zba1}).\end{proof}

\begin{remark} \label{r}
In this remark  we assume that the   KLS conjecture holds  and show that it implies  (\ref{bba1}).  The difference of this case from the proof of
Proposition \ref{Step0} is to improve  (\ref{at}) and (\ref{at2}) from the  time of order  $p^{-2}$  to the time   of order  $p^{-1}$, respectively.
The isoperimetric  constant bound  allows  us to  improve the estimate (2.23) in \cite{G24} so that  (\ref{ftf}) in Lemma \ref{Step4} above can be replaced  by
 \begin{align}
 \frac{d}{dt}  \mathbb{E} F_t\leq \big(500t^{-1}+K
 D_0^2  \big) \mathbb{E} F_t,\ \ \ \ \ t\in [t_0,1],
 \end{align}
where $K>0$ is a universal constant.  Then we can  replace  the estimate (\ref{at}) by
 \begin{align}
    \mathbb{E} F_{(9K)^{-1}p^{-1}}  \leq   & e^{-500 \log t_1^* + p/9 }   \mathbb{E} F_{t_1^*} \nonumber ,
   \end{align}
and also replace  the estimate (\ref{at2}) by
 \begin{align}
   \mathbb{E} F_{k_0,(9K)^{-1}p^{-1}}\leq
    e^{-500 \log t_1^* + p/9 }   \mathbb{E} F_{k_0,t_{k_0}} .\nonumber
   \end{align}
With the two estimates above, we may prove (\ref{bba1}) according to the procedure in the proof of Proposition   \ref{Step0} above.
Notice that the proof of   Proposition \ref{Step0}    depends on (\ref{kla}). This prevents  us from proving   (\ref{bba1}) in this paper. \medskip

\end{remark}

\begin{corollary}\label{Step000}
 There exists a universal constant $C>0$ such that for any isotropic  log-concave probability measure $\mu$ on   $\mathbb{R}^n$,
\begin{align}
 \mathbb{E}\mathrm{Tr}( A_t^p)
\leq (Cp)^{2p}n ,\ \ \ \ \ t\geq 0,\ p\geq 1.
\end{align}
\end{corollary} \noindent\textbf{Proof}\
Let $t=a/p^2$ for some $a\in (0,1)$. By Proposition  \ref{Step0} and (\ref{Lic}), we have for some   universal constant $c\in (0,1/15)$
\begin{align}
 \mathbb{E}\mathrm{Tr}( A_t^p)
 = & \sum_{ \lambda_{i;t}\leq \frac{8}{3}}  \mathbb{E}  \lambda_{i;t}^p+
 \sum_{ \lambda_{i;t}> \frac{8}{3}}  \mathbb{E}  \lambda_{i;t}^p
 \nonumber\\
 \leq &   3^pn+  e^{-ct^{-1/2}}t^{-p} n \nonumber\\
 \leq &   (3p)^pn+  e^{-(ca^{-1/2}+\log a)p}p^{2p} n.\nonumber
\end{align}
Since  $c\in (0,1/15)$, we can check that $ca^{-1/2}+\log a>0$ when $a<c^{4}$.
Therefore  we further have for $a <c^{4}$
\begin{align}
 \mathbb{E}\mathrm{Tr}( A_t^p)
  \leq &  (6p)^{2p}n,\ \ \ t=a/p^2 .\nonumber
\end{align}
When $a\geq c^4$, we have by  (\ref{Lic}) and $t\geq c^4p^{-2}$
\begin{align}
 \mathbb{E}\mathrm{Tr}( A_t^p)
\leq  (c^{-2}p)^{2p}n .\nonumber
\end{align}
Combing the two cases above, we finish the proof.
\qed

\section{   {Appendix} }

Let $\mu^*$  be an isotropic log-concave probability measure on $\mathbb{R}^k$ for some  $k\geq 2$ and denote by $(A^*_t)$ the corresponding covariance process of the stochastic localization.
 Let $0\leq \lambda_{1;t}^* \leq  \cdots \leq  \lambda_{k;t}^*$ be the eigenvalues of the covariance matrix $A^*_t$, repeated according to their multiplicity.
From the proof of the  Proposition 64 in \cite{KL24},    we can choose  $\mu^*$ such that, for  any $t\geq \frac{1}{\log k}$,
$\lambda_{k;t}^*\geq \frac{c}{t}$ with probability  no less than
$1-e^{-1/2}$. Here $c>0$ is some universal constant.

Assume that $n/k\geq 1$ is an integer and denote by  $A_t$
  the  covariance process of the stochastic localization for the probability  measure $\mu:=(\mu^*)^{\otimes \frac{n}{k}}$.
Denote by  $0\leq \lambda_{1;t} \leq  \cdots \leq  \lambda_{n;t}$   the eigenvalues of the covariance matrix $A_t$.
Let   $t = \frac{1}{\log k}$ and  $r=\frac{c}{t}$. Then we have
\begin{align}
\mathbb{E}\sum_{i=1, \lambda_{i;t}\geq r}^n\lambda_{i;t}^r\geq  (1-e^{-1/2}) \frac{n}{k}r^r=(1-e^{-1/2})e^{-r/c}r^rn.\nonumber
\end{align}
This shows that the conjecture in  (\ref{pm}) is sharp within a constant if it holds.

Next we show that it is   hard to give  the     estimate of $\mathbb{E} \mathrm{Tr}(A_t^p)$   in (\ref{pm})  by a modification of the method in \cite{G24}   for the case $p=2$.

First we prepare a special $p$.  By the calculation  above,  for $t = \frac{1}{\log k}$ and  $r=\frac{c}{t}$,  we   also have
\begin{align}
\mathbb{E}\sum_{i=1, \lambda_{i;t}\geq r}^n \lambda_{i;t}^q\geq  (1-e^{-1/2}) e^{-r/c} r^q n,\ \ \ \ \ q>0.\nonumber
\end{align}
This shows that  for   $t=\frac{c}{p^{3/2}}=\frac{1}{\log k}$
\begin{align}
\mathbb{E}\sum_{i=1, \lambda_{i;t}\geq p}^n \lambda_{i;t}^p\geq  (1-e^{-1/2}) e^{-p^{3/2}/c}  (p^{3/2} )^pn, \label{fa}
\end{align}
where    $p>1$ is  large enough such that  $ \frac{c}{p^{3/2}}<1/p$.

 We assume that $h$ is a    positive  test function such that
\begin{align}
h(u)=u^p,\ \ \ \ u\geq p. \label{hr}
\end{align}
Then we have by (\ref{fa}) and (\ref{hr})
\begin{align}
\mathbb{E}\sum_{i=1}^n h(\lambda_{i;t})\geq    (1-e^{-1/2}) e^{-p^{3/2}/c}  (p^{3/2} )^pn,\ \ \ \ \  t=\frac{c}{p^{3/2}}. \label{hr1}
\end{align}

Now we start at time $t=\frac{c}{p^{3/2}}<1/p$ to estimate $\mathbb{E}\sum_{i=1}^n h(  \lambda_{i;1/p}^p)$ in the example above. The best differential estimate in
 the proof of Lemma 2.4 in \cite{G24}    should be roughly
\begin{align}
 \frac{d}{ds}  \mathbb{E}\sum_{i=1}^n h(  \lambda_{i;s}^p) \leq & Cp(p-1)s^{-1}
  \mathbb{E}\sum_{i=1}^n h(  \lambda_{i;s}^p) ,\ \ \ s>0,\label{hr2}\end{align}
  which does not depend on the other  information  of  the test function except (\ref{hr}).
  Applying (\ref{hr1}) and (\ref{hr2}), we can get no more than
  \begin{align}
\mathbb{E}\sum_{i=1}^n h(\lambda_{i;1/ p})\leq   (1-e^{-1/2})e^{-p^{3/2}/c}  (p^{3/2} )^p e^{Cp(p-1)(\frac{1}{2}\log p-\log c)}n  .\nonumber
\end{align}
This is however   far from the estimate in the  conjecture.

\end{document}